\documentclass{amsart}
\usepackage{amssymb}
\def\C{{\mathbb C}}

\def\N{{\mathbb N}}
\def\R{{\mathbb R}}

\def\abs#1{\left|#1\right|}
\def\Alg{{\mathcal A}}
\def\Cnt#1{{\mathcal C}^{#1}}
\def\co#1{{#1}^c}
\def\comp{\circ}
\def\csub{\Subset}

\def\eps{\varepsilon}

\def\Gen{{\mathcal G}}
\def\Genld{\Gen_{\mathrm{ld}}}
\def\GenC{\widetilde\C}

\def\GenR{\widetilde\R}
\DeclareMathOperator\id{id}

\def\idealproper{\lhd}
\def\implies{\Rightarrow}

\def\interior#1{{#1}^\circ}
\def\inv#1{{#1}^{-1}}

\DeclareMathOperator\Ker{Ker}
\def\Mod{{\mathcal E}_M}

\def\Null{\mathcal N}
\def\Ord{\mathrm O}
\def\restr#1#2{{#1}_{|{#2}}}

\newtheorem{thm}{Theorem}[section]
\newtheorem{prop}[thm]{Proposition}
\newtheorem{lemma}[thm]{Lemma}
\newtheorem{cor}[thm]{Corollary}

\theoremstyle{definition}
\newtheorem*{df}{Definition}

\theoremstyle{remark}
\newtheorem*{rem}{Remark}

\begin{document}
\title{Isomorphisms of algebras of Colombeau generalized functions}

\author{Hans Vernaeve}
\address{Institut f\"ur Grundlagen der Bauingenieurwissenschaften\\
Technikerstra\ss e 13\\
A-6020 Innsbruck}
\thanks{Supported by research grants M949 and Y237 of the Austrian Science Foundation (FWF)}


\date{}

\begin{abstract}
We show that for smooth manifolds $X$ and $Y$, any isomorphism between the special algebra of Colombeau generalized functions on $X$, resp.\ $Y$ is given by composition with a unique Colombeau generalized function from $Y$ to $X$. We also identify the multiplicative linear functionals from the special algebra of Colombeau generalized functions on $X$ to the ring of Colombeau generalized numbers. Up to multiplication with an idempotent generalized number, they are given by an evaluation map at a compactly supported generalized point on $X$.
\end{abstract}

\maketitle

\section{Introduction}
It is a famous theorem in commutative Banach algebra theory that the isomorphisms between algebras of ($\C$-valued) continuous functions on compact Hausdorff topological spaces $X$, resp.~$Y$ are given by composition with a unique homeomorphism from $Y$ to $X$ \cite{Kadison}. When $X$, $Y$ are smooth Hausdorff manifolds, isomorphisms between algebras of smooth functions on them are similarly given by composition with a unique diffeomorphism from $Y$ to $X$. This result was only recently established in its full generality \cite{Grabowksi, Mrcun}.

A natural extension of these theorems is to look whether a similar theorem can hold in algebras of generalized functions, e.g., containing the space of distributions. Nonlinear generalized functions in the sense of J.F.~Colombeau \cite{Colombeau1984, DHPV2004} were introduced as a tool for studying nonlinear partial differential equations. They are an extension of the theory of distributions providing maximal consistency with respect to classical algebraic operations \cite{GKOS} in view of L.~Schwartz's impossibility result \cite{Schwartz}. Under the influence of applications of a more geometric nature (e.g.\ in Lie group analysis of differential equations and in general relativity), a geometric theory of Colombeau generalized functions arose \cite{GKOS, GKSV2002, KSV2003}. In particular, for $X$, $Y$ smooth paracompact Hausdorff manifolds, Colombeau generalized functions from $X$ to $Y$ can be defined. Recently, a definition of distributions from $X$ to $Y$ was proposed as a quotient of a subspace of the space $\Gen[X,Y]$ of so-called special Colombeau generalized functions from $X$ to $Y$ \cite{KSV2006}.

Denoting the algebra of (complex-valued) Colombeau generalized functions on a smooth paracompact Hausdorff manifold $X$ (resp.\ $Y$) by $\Gen(X)$ (resp.\ $\Gen(Y)$), we show more generally (theorem~\ref{main result} and its corollary) that algebra homomorphisms $\Gen(X)\to\Gen(Y)$ are characterized as compositions with locally defined Colombeau generalized functions from $X$ to $Y$, up to multiplication with an idempotent element of $\Gen(Y)$ (which is necessarily locally constant on $Y$). When the homomorphism is an isomorphism, the idempotent element necessarily equals $1$ and the generalized function from $X$ into $Y$ is uniquely determined.

Our technique is based on a characterization of the multiplicative $\GenC$-linear functionals on $\Gen(X)$, where $\GenC$ denotes the ring of Colombeau generalized complex numbers. Up to multiplication with an idempotent element of $\GenC$, these functionals coincide with the evaluation maps at generalized points \cite[\S 3.2]{GKOS} in $\Gen(X)$.

\section{Preliminaries}
The ring $\GenC$ of (complex) Colombeau generalized numbers is defined as $\Mod/\Null$, where
\begin{align*}
\Mod =\,&\{(z_\eps)_\eps \in \C^{(0,1)}: (\exists b\in\R) (\abs{z_\eps}=\Ord(\eps^b),\text{ as }\eps\to 0 )\}\\
\Null =\,&\{(z_\eps)_\eps \in\Mod: (\forall b\in\R) (\abs{z_\eps}=\Ord(\eps^b),\text{ as }\eps\to 0 )\}.
\end{align*}
Colombeau generalized numbers arise naturally as evaluations of a Colombeau generalized function at a point in its domain. The subring of $\GenC$ consisting of those elements that have a net of real numbers as a representative, is denoted by $\GenR$. Nets in $\Mod$ are called moderate, nets in $\Null$ negligible. $\GenC$ is a complete topological ring with zero divisors; the associated topology is called the sharp topology \cite{AJ2001}.

Let $S\subseteq(0,1)$. By $e_S\in\GenR$ we denote the element which has the characteristic function on $S$ as a representative. Then every idempotent element in $\GenC$ is of the form $e_S$ \cite{AJOS2006}. Further algebraic properties of $\GenC$ are described in \cite{AJ2001, AJOS2006}.

By a smooth manifold, we will mean a second countable Hausdorff $\Cnt{\infty}$ manifold of finite dimension (without boundary).

Let $X$ be a smooth manifold. By $K\csub X$, we denote a compact subset $K$ of $X$.
Let $\xi\in\mathcal X(X)$ denote a vector field on $X$ and $L_\xi$ its Lie derivative.
Then the (so-called special) algebra $\Gen(X)$ of Colombeau generalized functions on $X$ is defined as $\Mod(X)/\Null(X)$, where
\begin{align*}
\Mod(X) =\,&\big\{(u_\eps)_\eps \in(\Cnt{\infty}(X))^{(0,1)}: (\forall
K\csub X) (\forall k\in\N) (\exists b\in\R)\\&(\forall \xi_1,\ldots, \xi_k\in {\mathcal X}(X))
\big(
\sup_{p\in K}\abs{L_{\xi_1}\cdots L_{\xi_k} u_\eps(p)}=\Ord(\eps^b),\text{ as }\eps\to 0
\big)\big\}\\
\Null(X) =\,&\big\{(u_\eps)_\eps \in\Mod(X): (\forall
K\csub X) (\forall b\in\R)
\big(
\sup_{p\in K}\abs{u_\eps(p)}=\Ord(\eps^b),\text{ as }\eps\to 0
\big)\big\}.
\end{align*}
See also \cite[\S 3.2]{GKOS} for several equivalent definitions.

A net $(p_\eps)_\eps\in X^{(0,1)}$ is called compactly supported \cite[\S 3.2]{GKOS} if there exists $K\csub X$ and $\eps_0>0$ such that $p_\eps\in K$, for $\eps<\eps_0$. Denoting by $d_h$ the Riemannian distance induced by a Riemannian metric $h$ on $X$, two nets $(p_\eps)_\eps$, $(q_\eps)_\eps$ are called equivalent if the net $(d_h(p_\eps, q_\eps))_\eps$ is negligible (this does not depend on the choice of $h$). The equivalence classes w.r.t.\ this relation are called compactly supported generalized points on $X$. The set of compactly supported generalized points on $X$ will be denoted by $\widetilde X_c$. If $u\in\Gen(X)$ and $p\in \widetilde X_c$, the point value $u(p)\in\GenC$ is the generalized number with representative $(u_\eps(p_\eps))_\eps$ (this does not depend on the representatives).

Let $X$, $Y$ be smooth manifolds. The space $\Gen[X,Y]$ of c-bounded Colombeau generalized functions from $X$ to $Y$ is similarly defined as a quotient of the set $\Mod[X,Y]$ of moderate, c-bounded nets of smooth maps $X\to Y$ (\cite[Def.~3.2.44]{GKOS}) by a certain equivalence relation $\sim$ (\cite[Def.~3.2.46]{GKOS}). (See also Appendix for the notion of c-boundedness.)

We will also use a slightly modified version of the space $\Gen[X,Y]$ where we do not require the nets to be globally defined. The space $\Genld[X,Y]$ of locally defined c-bounded Colombeau generalized functions $X\to Y$ is the set of all nets $(u_\eps)_\eps$ of smooth maps defined on $X_\eps\subseteq X\to Y$ with the property that $(\forall K\csub X)(\exists \eps_0>0)(\forall\eps<\eps_0)(K\subseteq X_\eps)$ and satisfying the c-boundedness and moderateness conditions for elements of $\Mod[X,Y]$, modulo the equivalence relation $\sim$ as it is defined on $\Mod[X,Y]$. 
By definition, $\Gen[X,Y]$ is a subset of $\Genld[X,Y]$.
\begin{rem}
Under mild topological restrictions on $X$, $\Genld[X,Y]=\Gen[X,Y]$. E.g., it is sufficient that
$(\forall K\csub X) (\exists f\in\Cnt{\infty}(X,X))(\overline{f(X)}\csub X \,\&\, \restr{f}{K}=\id_K)$. This appears to be fulfilled in almost all practical cases.
\end{rem}

\section{Surjectivity of multiplicative $\GenC$-linear maps}
Throughout this paper, $\Alg$, $\Alg_1$, $\Alg_2$,\dots\ are commutative $\GenC$-algebras with $1$.
By a linear map $\Alg_1\to\Alg_2$, a $\GenC$-linear map is meant.
In particular, a multiplicative linear functional on $\Alg$ is meant to be a multiplicative $\GenC$-linear map $\Alg\to\GenC$.

\begin{lemma}\label{surjective and 1}
\hfil
\begin{enumerate}
\item If a multiplicative linear map $\phi$: $\Alg_1\to\Alg_2$ is surjective, then $\phi(1)=1$.
\item A multiplicative linear functional $m$ on $\Alg$ is surjective iff $m(1)=1$.
\end{enumerate}
\end{lemma}
\begin{proof}
Elementary.
\end{proof}

\begin{prop}\label{surjective-decomposition}
Suppose that there exists a multiplicative linear map $\phi_0$: $\Alg_1\to\Alg_2$ with $\phi_0(1)=1$.
Let $\phi$ be any multiplicative linear map $\Alg_1\to\Alg_2$. Then there exists a multiplicative linear map $\psi$: $\Alg_1\to\Alg_2$ with $\psi(1)=1$ such that $\phi = \phi(1)\cdot\psi$.\\
If $\Alg$ is a topological algebra and $\phi_0$, $\phi$ are continuous, then $\psi$ is also continuous.
\end{prop}
\begin{proof}
Let $\psi = \phi + (1- \phi(1)) \phi_0$.
\end{proof}
E.g., if $\Alg_1=\Gen(X)$, $X$ a manifold, then for any $p\in X$, $\delta_p$: $\Alg_1\to\Alg_2$: $\delta_p(u)=u(p)1$ is a multiplicative linear map $\Alg_1\to\Alg_2$ with $\delta_p(1)=1$.
In particular, the study of multiplicative linear functionals on $\Gen(X)$ is reduced to the surjective ones.

\section{Multiplicative $\GenC$-linear functionals on $\Gen(X)$}
For a (non-zero) multiplicative $\C$-linear functional $m$ on a $\C$-algebra $A$, $A/\Ker m$ $\cong \C$ is a field, so $\Ker m$ is a maximal ideal. If $A$ is a Banach algebra, the converse also holds: for a maximal ideal $M\idealproper A$, $A/M\cong \C$ by the Gelfand-Mazur theorem \cite[3.2.4]{Kadison}, and the canonical surjection $A\to A/M$ determines a multiplicative $\C$-linear functional. Since $\GenC$ is not a field, the kernel of a multiplicative $\GenC$-linear functional on a $\GenC$-algebra $\Alg$ will not be a maximal ideal. This motivates the following definition.

\begin{df}
An ideal $I\idealproper\Alg$ is maximal with respect to the property $I\cap\GenC 1=\{0\}$ iff $J\idealproper\Alg$, $I\subseteq J$ and $J\cap\GenC 1 =\{0\}$ imply that $I=J$.
\end{df}
It is easy to see that for a surjective multiplicative $\GenC$-linear functional $m$ on $\Alg$, $\Ker m$ is an ideal maximal with respect to $\Ker m\cap\GenC 1 = \{0\}$.

\begin{df}
Let $u\in\Alg$ and $S\subseteq (0,1)$ with $0\in\overline S$. Then $u$ is called invertible w.r.t.\ $S$ iff
there exists $v\in\Alg$ such that $uv=e_S$.
\end{df}

\begin{lemma}
Let $I\idealproper\Alg$. The following are equivalent:
\begin{enumerate}
\item $I\cap\GenC 1=\{0\}$
\item for each $S\subseteq(0,1)$ with $0\in\overline S$, if $u\in\Alg$ and $u$ is invertible with respect to $S$, then $u\notin I$.
\end{enumerate}
\end{lemma}
\begin{proof}
Let $u$ be invertible w.r.t.\ $S$. Should $u\in I$, then also $0\ne e_S\in I$, so $I\cap\GenC 1=\{0\}$.
\end{proof}
We denote the complement of $S\subseteq (0,1)$ by $\co S$.
\begin{lemma}\label{invertible wrt S}
Let $u\in\Alg$ and $S\subseteq (0,1)$ with $0\in\overline S$. Then $u$ is invertible w.r.t.\ $S$ iff $u e_S + e_{\co S}$ is invertible.
\end{lemma}
\begin{proof}
If $uv=1$, for some $v\in \Alg$, then $(ue_S+e_{\co S})(ve_S + e_{\co S})= uve_S + e_{\co S}=1$.\\
Conversely, if $(ue_S+e_{\co S})v=1$, for some $v\in\Alg$, then multiplying by $e_S$ shows that $u (v e_S) = e_S$.
\end{proof}
\begin{cor}\label{cor invertible wrt S}
Let $X$ be a smooth submanifold of $\R^d$.
Let $S\subseteq (0,1)$ with $0\in\overline S$. Let $u\in\Gen(X)$. Then the following are equivalent:
\begin{enumerate}
\item $u$ is invertible w.r.t.\ $S$ (as an element of $\Gen(X)$)
\item $u(\tilde x)$ is invertible w.r.t.\ $S$ (as an element of $\GenC$), for each $\tilde x\in\widetilde X_c$.
\end{enumerate}
\end{cor}
\begin{proof}
This is a combination of the previous lemma with proposition~\ref{pointwise invertibility}.
\end{proof}

\begin{prop}\label{key lemma}
Let $X$ be a smooth submanifold of $\R^d$. Let $I\idealproper\Gen(X)$. If $(\forall p\in\widetilde X_c)(\exists u_p\in I)(u_p(p)\ne 0)$, then $I$ is not the kernel of a surjective multiplicative linear functional on $\Gen(X)$.
\end{prop}
\begin{proof}
Suppose that $I$ is the kernel of a surjective multiplicative linear functional. Then $I\cap\GenC 1=\{0\}$ and $I+\GenC 1 =\Gen(X)$, so each of the functions $x_i\in I+\GenC 1$ ($i\in\{1,\dots,d\}$), i.e., for each $i$, there exists $\lambda_i\in\GenC$ such that $x_i-\lambda_i 1\in I$. Write $\lambda=(\lambda_1,\dots,\lambda_d)\in\GenC^d$ and consider $\abs{x-\lambda}^2=\sum_i (x_i-\overline\lambda_i 1)(x_i-\lambda_i 1)\in I$.\\
We distinguish 3 cases.\\
(1) $\lambda\in\widetilde X_c$. Notice that by corollary~\ref{cor equivalent metrics}, this property is well-defined. Then also $\abs{x-\lambda}^2 + \abs{u_\lambda}^2\in I$. As $u_\lambda(\lambda)\ne 0$, there exists $S\subseteq (0,1)$, $0\in \overline S$, such that $u_\lambda(\lambda)\in\GenC$ is invertible w.r.t.\ $S$. 
Let $\tilde x\in\widetilde X_c$ with representative $(x_\eps)_\eps$. By proposition~\ref{sharp topology}, there exist $m$, $k\in\N$ such that
\[(\exists \eps_0>0)(\forall\eps\in S\cap(0,\eps_0))
(\abs{x_\eps-\lambda_\eps}\le\eps^m\implies\abs{u_{\lambda,\eps}(x_\eps)}\ge \eps^k).\]
So $(\exists \eps_0>0)$ $(\forall\eps\in S\cap(0,\eps_0))$ 
$
\bigg(\abs{x_\eps-\lambda_\eps}^2 +\abs{u_{\lambda,\eps}(x_\eps)}^2\ge
\begin{cases}
\eps^{2k},& \abs{x_\eps-\lambda_\eps}\le\eps^m\\
\eps^{2m},& \abs{x_\eps-\lambda_\eps}\ge\eps^m
\end{cases}
\bigg),
$
and we conclude by corollary~\ref{cor invertible wrt S} that $\abs{x-\lambda}^2 + \abs{u_\lambda}^2\in I$ is invertible w.r.t.\ $S$, a contradiction.\\
(2) $\lambda\in\widetilde X\setminus\widetilde X_c$, where $\widetilde X=\{\tilde x\in\GenR^d: (\exists $ repr.\ $(x_\eps)_\eps$ of $\tilde x)(\forall\eps)(x_\eps\in X)\}$. Let $(K_n)_{n\in\N}$ be a compact exhaustion of $X$ with $K_n\subseteq\interior{(K_{n+1})}$, $\forall n\in\N$ (where the interior is taken in the relative topology on $X$). Consider a representative $(\lambda_\eps)_\eps$ of $\lambda$ such that $\lambda_\eps\in X$, $\forall \eps$. As $\lambda\notin \widetilde X_c$, there exists a decreasing sequence $(\eps_n)_{n\in\N}$ with $\eps_n\to 0$ such that $\lambda_{\eps_n}\in X\setminus K_n$ for each $n$. As the Euclidean distance $d(X\setminus K_n, K_{n-1})>0$ for each $n$, $v(x)=\abs{x-\lambda_\eps}^2\in I$ evaluated in any compactly supported point of $X$ is invertible w.r.t.\ $S=\{\eps_n: n\in\N\}$, a contradiction.\\
(3) If $\lambda\in\GenC^d\setminus \widetilde X$, then for any representative $(\lambda_\eps)_\eps$ of $\lambda$, $(d(\lambda_\eps, X)_\eps)_\eps$ is not a negligible net. This means that there exists $S\subseteq (0,1)$ with $0\in\overline S$ and $m\in\N$ such that $d(\lambda_\eps, X)\ge \eps^m$, for each  $\eps\in S$. This also means that $v(x)=\abs{x-\lambda_\eps}^2\in I$ evaluated in any compactly supported point of $X$ is invertible w.r.t.\ $S$, a contradiction.
\end{proof}


\begin{thm}\label{mlf characterization}
Let $X$ be a smooth manifold.
\begin{enumerate}
\item The surjective multiplicative linear functionals on $\Gen(X)$ are
\[\delta_p: \Gen(X)\to\GenC: \delta_p(u) = u(p),\]
where $p\in\widetilde X_c$.
\item The multiplicative linear functionals on $\Gen(X)$ are
\[e\delta_p: \Gen(X)\to\GenC: e\delta_p(u)=e u(p),\]
where $p\in\widetilde X_c$ and $e\in\GenR$ idempotent.
\end{enumerate}
\end{thm}
\begin{proof}
First, let $X$ be a smooth submanifold of $\R^d$.\\
(1) Let $m$ be a surjective multiplicative linear functional on $\Gen(X)$. Then by the lemma, there exists $p\in\widetilde X_c$ such that $u(p)=0$, $\forall u\in \Ker m$. I.e., $\Ker m\subseteq \Ker \delta_p$. But $\Ker m$ is maximal w.r.t.\ $\Ker m\cap\C 1=\{0\}$ and $\Ker \delta_p\cap \C 1 =\{0\}$, so $\Ker m=\Ker\delta_p$. So for each $u\in\Gen(X)$, as $u-u(p)\in\Ker m=\Ker\delta_p$, $m(u)=m(u-u(p)+u(p))=m(u(p))=u(p)$, so $m=\delta_p$.\\
(2) This follows from part~1 and proposition~\ref{surjective-decomposition}.\\
Now let $X$ be any smooth manifold. It follows from Whitney's embedding theorem \cite{Hirsch} that there exists a smooth embedding $f$: $X\to\R^d$, for some $d\in\N$.
Let $m$: $\Gen(X)\to\GenC$ be a surjective multiplicative linear functional. For $u\in\Gen(f(X))$, $u\comp f\in\Gen(X)$ (corollary \ref{cor c-bounded}). Then $\mu$: $\Gen(f(X))\to\GenC$: $\mu(u)=m(u\comp f)$ is a surjective multiplicative linear functional, so there exists $p\in\widetilde{f(X)}_c$ such that $\mu(u)=u(p)$, $\forall u\in\Gen(f(X))$. For each $v\in\Gen(X)$, $v\comp \inv{f}\in\Gen(f(X))$, so $m(v)=\mu(v\comp\inv{f})=v(\inv{f}(p))$, where $\inv{f}(p)\in\widetilde X_c$ \cite[3.2.55]{GKOS}.
\end{proof}

\section{Algebra homomorphisms $\Gen(X)\to\Gen(Y)$}
\begin{thm}\label{main result}
Let $X\subseteq\R^{d_1}$, $Y\subseteq\R^{d_2}$ be smooth submanifolds.
\begin{enumerate}
\item Let $\phi$: $\Gen(X)\to\Gen(Y)$ be a morphism of algebras (i.e., a multiplicative $\GenC$-linear map).
Then there exists $f\in(\Gen(Y))^{d_1}$, c-bounded into $X$ and $e\in\Gen(Y)$ idempotent such that
\[\phi(u)=e\cdot (u\comp f), \quad\forall u\in\Gen(X).\]
If $\phi(1)=1$, then $e=1$ 
and $f$ is uniquely determined.
\item If $\phi$: $\Gen(X)\to\Gen(Y)$ is an isomorphism of algebras (i.e., additionally, $\phi$ is bijective), then the map $f$ has an inverse $\inv f \in(\Gen(X))^{d_2}$, c-bounded into $Y$ such that $\inv\phi$ is given by composition with $\inv f$. As a map $\widetilde Y_c\to\widetilde X_c$, $f$ is bijective. In this case, $\dim X = \dim Y$.
\end{enumerate}
\end{thm}
\begin{proof}
(1) First, let $\phi(1)=1$. Let $\tilde x\in\widetilde Y_c$ arbitrary. Then the map $\delta_{\tilde x}\comp \phi$ is a multiplicative linear functional on $\Gen(X)$. It is also surjective, as $\delta_{\tilde x}(\phi(1))=1$.
So by theorem~\ref{mlf characterization}, there exists $f(\tilde x)\in\widetilde X_c$ such that $\delta_{\tilde x}\comp\phi=\delta_{f(\tilde x)}$. So
\begin{equation}\label{concrete delta equality}
(\forall u\in\Gen(X)) (\forall \tilde x\in\widetilde Y_c) ((\phi(u))(\tilde x)= u(f(\tilde x))).
\end{equation}
In particular, for $u_i(x)=x_i\in\Gen(X)$, $i=1,\dots, d_1$, we see that
\begin{equation}\label{function f}
(\phi(u_1),\dots, \phi(u_{d_1}))\in(\Gen(Y))^{d_1}
\end{equation}
is the unique generalized function which coincides with $f$ when evaluated at generalized points in $\widetilde Y_c$ (because an element of $\Gen(Y)$ is completely determined by its values in $\widetilde Y_c$ \cite[Thm.~3.2.8]{GKOS}). With a slight abuse of notation, we will therefore also denote it by $f$. By proposition~\ref{c-bounded}, $f$ is c-bounded into $X$. So by proposition~\ref{c-bounded composition}, for each $u\in\Gen(X)$, the componentwise composition $u\comp f$ defines an element of $\Gen(Y)$. By eqn.~(\ref{concrete delta equality}), it coincides with $\phi(u)$ on each compactly supported point in $\widetilde Y_c$, so $u\comp f = \phi(u)$ in $\Gen(Y)$. Clearly, $f$ is completely determined by $f_i=u_i\comp f=\phi(u_i)$ ($i=1$, \dots, $d_1$).\\
For general $\phi$, this follows by proposition~\ref{surjective-decomposition} and the fact that $\phi(1)$ is idempotent.\\
(2) Applying part~1 on $\inv\phi$, we find $g\in(\Gen(X))^{d_2}$, c-bounded into $B$ such that $\inv\phi$ is given by composition with $g$. To see that $g=\inv f$, we show that $f\comp g=\id_{\Gen(X)}\in(\Gen(X))^{d_1}$, where $\id_{\Gen(X)}$ is the generalized function with representative $(\id_X)_\eps$.\\
By eqn.~(\ref{function f}) and because $\inv\phi$ is given by composition with $g$,
\begin{multline*}
f \comp g = (f_1\comp g,\dots, f_{d_1}\comp g)
=(\inv\phi(f_1),\dots, \inv\phi(f_{d_1}))\\
= (\inv\phi(\phi(u_1)),\dots, \inv\phi(\phi(u_{d_1})))
=(u_1,\dots,u_{d_1}).
\end{multline*}
Similarly, $g\comp f=\id_{\Gen(Y)}\in(\Gen(Y))^{d_2}$. From these equalities, it follows also that $\inv f$ is the inverse of $f$ as pointwise maps on compactly supported generalized points.\\
Suppose that $m=\dim X < \dim Y =n$.
Let $y\in Y$ and  let $W$ be a geodesically convex neighbourhood of $y$ in $Y$ with $\overline W\csub U$ for some chart $(U,\phi)$ of $Y$. By the c-boundedness of $f$, $f_\eps(W)\subseteq L\csub X$ for sufficiently small $\eps$. We can cover $L$ by open charts of $X$. By compactness, the open cover has a Lebesgue number $\delta>0$.
There exists $M\in\N$ such that for sufficiently small $\eps$ and for each $y'\in W$,
$\abs{f_\eps(y)-f_\eps(y')}\le \frac{\eps^{-M}}{2}\abs{\phi(y)-\phi(y')}$.
So for sufficiently small $\eps$, $f_\eps\comp \inv\phi(B(\phi(y),{\delta}\eps^M))\subseteq V_\eps$ for some chart $(V_\eps,\psi_\eps)$ of $X$. Applying the Borsuk-Ulam theorem \cite[\S 16.5]{Zeidler} on $\psi_\eps\comp f_\eps \comp \inv\phi$: $B(\phi(y), {\delta}\eps^M)\to \R^m$, we obtain $p_\eps$, $p'_\eps$ $\in W$ with $\abs{\phi(p_\eps)-\phi(p'_\eps)}\ge{\delta}\eps^M$ such that $f_\eps(p_\eps)=f_\eps(p'_\eps)$ by the injectivity of $\psi_\eps$. By \cite[Lemma~3.2.6]{GKOS}, $(p_\eps)_\eps$, $(p'_\eps)_\eps$ represent different elements of $\widetilde Y_c$ with equal images under $f$, a contradiction. Similarly, $\dim X\le\dim Y$.
\end{proof}
\begin{cor}
Let $X$, $Y$ be smooth manifolds.
\begin{enumerate}
\item Let $\phi$: $\Gen(X)\to\Gen(Y)$ be a morphism of algebras (i.e., a multiplicative $\GenC$-linear map).
Then there exists $f\in\Genld[Y,X]$ and $e\in\Gen(Y)$ idempotent such that
\[\phi(u)=e\cdot (u\comp f), \quad\forall u\in\Gen(X).\]
If $\phi(1)=1$, then $e=1$ 
and $f$ is uniquely determined.
\item If $\phi$: $\Gen(X)\to\Gen(Y)$ is an isomorphism of algebras (i.e., additionally, $\phi$ is bijective), then the map $f$ has an inverse $\inv f \in\Genld[X,Y]$ such that $\inv\phi$ is given by composition with $\inv f$. As a map $\widetilde X_c\to\widetilde Y_c$, $f$ is bijective. In this case, $\dim X = \dim Y$.
\end{enumerate}
\end{cor}
\begin{proof}
It follows from Whitney's embedding theorem \cite{Hirsch} that there exist smooth embeddings $\iota_1$: $X\to\R^{d_1}$ and $\iota_2$: $Y\to\R^{d_2}$, for some $d_1$, $d_2\in\N$.\\
(1) Let $\phi$: $\Gen(X)\to\Gen(Y)$ be a  multiplicative $\GenC$-linear map with $\phi(1)=1$. Then $\widetilde\phi$: $\Gen(\iota_1(X))\to\Gen(\iota_2(Y))$: $\widetilde\phi(u)=\phi(u\comp\iota_1)\comp\inv\iota_2$ is a multiplicative $\GenC$-linear map with $\widetilde\phi(1)=1$.
By the previous theorem and by corollary~\ref{cor c-bounded}, there exists $\widetilde f\in\Genld[\iota_2(Y),\iota_1(X)]$ such that $\widetilde\phi$ is given by composition with $\widetilde f$. So for each $u\in\Gen(X)$, $\phi(u)=\widetilde\phi(u\comp\inv\iota_1)\comp\iota_2=u\comp(\inv\iota_1\comp \widetilde f\comp\iota_2)$. By the analogue of \cite[Cor.~3.2.59]{GKOS} for $\Genld[X,Y]$, $f=\inv\iota_1\comp \widetilde f\comp\iota_2\in \Genld[Y,X]$. Unicity of $f$ follows from unicity of $\widetilde f$.\\
The result for general $\phi$ follows again from proposition~\ref{surjective-decomposition}.\\
(2) We similarly find $\widetilde g\in\Genld[\iota_1(X),\iota_2(Y)]$ with $g=\inv\iota_2\comp \widetilde g\comp\iota_1\in\Genld[X,Y]$. By the previous theorem, $\widetilde f\comp\widetilde g$ is the identity in $\Genld[\iota_1(X),\iota_1(X)]$. So $f\comp g=\inv\iota_1\comp\id_{\Genld[\iota_1(X),\iota_1(X)]}\comp\iota_1=\id_{\Genld[X,X]}$, and similarly, $g\comp f=\id_{\Genld[Y,Y]}$. It follows again that $g=\inv f$ as pointwise maps on compactly supported generalized points.
\end{proof}

Concerning idempotent elements in $\Gen(X)$, we can be more explicit:
\begin{prop}
Let $X$ be a smooth manifold. Let $e\in\Gen(X)$ be an idempotent element (i.e., $e^2=e$ holds in $\Gen(X)$). Then on every connected component of $X$, $e$ is an idempotent constant.
\end{prop}
\begin{proof}
If $X$ is an open subset of $\R^d$, this is proven in \cite{AJOS2006}. Let $X$ be an arbitrary manifold. Consider a chart $(V,\psi)$ of $X$ and $x\in V$. Then the local representation $e\comp\inv\psi\in\Gen(\psi(V))$ is an idempotent, and therefore equal to some constant $c\in\GenC$ in a connected, open neighbourhood $W$ of $\psi(x)$. So $e=c$ in the open neighbourhood $\inv\psi(W)$ of $x$. Therefore, for every $c\in\GenC$, $\{x\in X: (\exists U$ open neighbourhood of $x)(\restr{e}{U}=c)\}$ is open and closed in $X$. Consequently, on every connected component $C$ of $X$, each $x\in C$ has an open neighbourhood $U$ such that $\restr{e}{U}=c$, for some constant $c\in\GenC$ independent of $x\in C$. The proposition follows by the fact that $\Gen(C)$ is a sheaf of differential algebras on $C$ (\cite[Prop.~3.2.3]{GKOS}).
\end{proof}

\appendix
\section{Colombeau generalized functions on a manifold embedded in $\R^d$}
In this appendix, we extend some results that are well-known in the special case where $X$ is an open subset of $\R^d$ to the case of a submanifold of $\R^d$.

\begin{lemma}\label{equivalent metrics}
Let $X$ be a connected smooth submanifold of $\R^d$.
Let $h$ be the Riemannian metric on $X$ induced by the Euclidean metric in $\R^d$. Let $K\csub X$.
Then there exists $C\in\R^+$ such that for each $p$, $q$ $\in K$, $\abs{p-q}\le d_h(p, q)\le C\abs{p-q}$.
\end{lemma}
\begin{proof}
$d_h(p,q)$ is the infimum of the distances between $p$, $q$ along paths on $X$, and therefore at least equal to the Euclidean distance between $p$ and $q$. For the other inequality, suppose first that $p$, $q$ lie in a sufficiently small 
neighbourhood 
of a given point $p_0\in K$. It is an exercise in elementary differential geometry that in this case, $d_h(p,q)\le C \abs{p-q}$ (with $C\to 1$ as $p$, $q$ $\to p_0$).\\
If the inequality would not hold globally on $K$, one could construct sequences $(p_m)_m$, $(q_m)_m$ of points in $K$ such that $d_h(p_m,q_m)\ge m \abs{p_m-q_m}$. Because $K$ is compact, there is a subsequence $(m_k)_k$ such that $p_{m_k}\to p\in K$, $q_{m_k}\to q\in K$. By continuity, $d_h(p,q)\ge m\abs{p-q}$, for each $m\in\N$, so $\abs{p-q}=0$ and $p=q$. This contradicts the inequality in an arbitrary small neighbourhood of $p$.
\end{proof}
\begin{cor}\label{cor equivalent metrics}
Let $X$ be a smooth submanifold of $\R^d$.
The compactly supported generalized points in $\widetilde X_c$ are in 1-1 correspondence with the elements of $\GenR^d$ which have a representative that consists of elements of $K$, for some $K\csub X$. More specifically, the injection is given by the (well-defined) map $\widetilde X_c\to \GenR^d$ which is the identity-map on representatives.
\end{cor}
\begin{proof}
By the fact that every $x\in X$ has a connected neighbourhood, $K\csub X$ is contained in a finite number of connected components of $X$.\\
Two compactly supported nets $(p_\eps)_\eps$, $(q_\eps)_\eps$ in $X^{(0,1)}$ represent the same generalized point in $X$ iff $d_h(p_\eps,q_\eps)=\Ord(\eps^m)$, $\forall m\in\N$. 
(By definition, this also implies that for a fixed sufficiently small $\eps$, $p_\eps$ and $q_\eps$ lie in the same connected component.)
By lemma \ref{equivalent metrics}, this is equivalent with $\abs{p_\eps-q_\eps}=\Ord(\eps^m)$, $\forall m\in\N$ (this also implies that for a fixed sufficiently small $\eps$, $p_\eps$ and $q_\eps$ lie in the same connected component, since any open cover of $K\csub X$, in particular one that consists of connected sets, has a Lebesgue number), i.e., they represent the same element in $\GenR^d$.
\end{proof}

\begin{prop}\label{pointwise invertibility}
Let $X$ be a smooth submanifold of $\R^d$.
Let $u\in\Gen(X)$. Then the following are equivalent:
\begin{enumerate}
\item $u$ is invertible (as an element of $\Gen(X)$)
\item $u(\tilde x)$ is invertible (as an element of $\GenC$), for each $\tilde x\in\widetilde X_c$.
\end{enumerate}
\end{prop}
\begin{proof}
$(1)\implies(2)$ is analogous to \cite[Thm.~1.2.5]{GKOS}.\\
$(2)\implies(1)$: to show that a global inverse exists, it is enough to show that there exists an inverse in each local representation (w.r.t.\ charts), and that the compatibility-conditions between them are satisfied \cite[Prop.~3.2.3]{GKOS}. As in \cite[Thm.~3.2.8]{GKOS}, part~(2) is also satisfied for each local representation. So by \cite[Thm.~1.2.5]{GKOS}, local inverses exist. The compatibility-conditions for $\inv{u}$ follow from the compatibility-conditions of $u$ and the fact that inverses in $\Gen(X)$ are unique (for any manifold $X$).
\end{proof}

\begin{prop}[Continuity in the sharp topology]\label{sharp topology}
Let $X$ be a smooth submanifold of $\R^d$. Let $u\in\Gen(X)$ and let $K\csub X$. Then for each $k\in\N$,
\begin{multline*}
(\exists m\in\N)(\exists \eps_0>0)(\forall\eps\le\eps_0)
(\forall x,y\in K)(\abs{x-y}\le\eps^m\implies\abs{u_{\eps}(x)-u_\eps(y)}\le \eps^k).
\end{multline*}
\end{prop}
\begin{proof}
If $X$ is an open subset of $\R^d$, see e.g.\ \cite[Prop.~3.1]{OPS2003}.\\
If $X$ is a smooth manifold of $\R^d$, cover $K$ by geodesically convex $W_\alpha$ with $\overline W_\alpha\csub V_\alpha$ for charts $(V_\alpha, \psi_\alpha)$ (as in \cite[Thm.~3.2.8]{GKOS}). By compactness, a finite number $W_1$, \dots, $W_M$ is sufficient. Call the corresponding charts $(V_1, \psi_1)$, \dots, $(V_M, \psi_M)$. By the existence of a Lebesgue number, we may suppose that $x$ and $y$ belong to the same $W_i$, if $\eps_0$ is chosen sufficiently small (and $m\ge 1$). So, let $k\in\N$. We apply the proposition to 
$u\comp\inv\psi_i\in\Gen(\psi_i(V_i))$, and we obtain $m_i\in\N$, $\eps_i> 0$ such that
\[
(\forall\eps\le\eps_i)
(\forall x,y\in W_i)(\abs{\psi_i(x)-\psi_i(y)}\le\eps^{m_i}\implies\abs{u_{\eps}(x)-u_\eps(y)}\le \eps^k).
\]
Further, by \cite[Lemma~3.2.6]{GKOS} and lemma~\ref{equivalent metrics} (as $W_i$ is connected),
$\abs{\psi_i(x)-\psi_i(y)}\le C d_h(x, y)\le C' \abs{x - y}$,
for some $C$, $C'$ $\in\R^+$ (independent of $x$, $y$ $\in W_i$). So, possibly after increasing $m_i$ and decreasing $\eps_i$,
\[
(\forall\eps\le\eps_i)
(\forall x,y\in W_i)(\abs{x-y}\le\eps^{m_i}\implies\abs{u_{\eps}(x)-u_\eps(y)}\le \eps^k).
\]
Choose $\eps_0\le \eps_1$, \dots, $\eps_0\le\eps_M$ and $m\ge m_1$, \dots, $m\ge m_M$. Then we obtain the statement of the proposition.
\end{proof}

Let $X\subseteq\R^{d_1}$, $Y\subseteq\R^{d_2}$ be smooth submanifolds. In analogy with the case where $X$, $Y$ are open subsets of $\R^{d_1}$, resp.~$\R^{d_2}$ (\cite[1.2.7]{GKOS}), $u\in (\Gen(X))^{d_2}$ is called c-bounded into $Y$ if there exists a representative $(u_\eps)_\eps$ of $u$ such that
\begin{equation}\label{eqn c-bounded}
(\forall K\csub X)(\exists K'\csub Y)(\exists \eps_0>0)(\forall \eps\le\eps_0)(u_\eps(K)\subseteq K').
\end{equation}

\begin{prop}\label{c-bounded composition}
Let $X\subseteq\R^{d_1}$, $Y\subseteq\R^{d_2}$ be smooth submanifolds. Let $u\in(\Gen(X))^{d_2}$ be c-bounded into $Y$ and $v\in\Gen(Y)$. Then the composition $v\comp u$ defined on representatives by means of
$(v\comp u)_\eps = (v_\eps\comp u_\eps)$
is a well-defined generalized function in $\Gen(X)$.
\end{prop}
\begin{proof}
Notice that the net $(v_\eps\comp u_\eps)_\eps$ is only locally defined; to find a globally defined representative, it can be multiplied by a net $(\chi_\eps)_\eps$ of smooth, compactly supported cut-off functions which is a representative of $1\in\Gen(X)$. Well-definedness 
follows as in \cite[Prop.~1.2.8]{GKOS}.
\end{proof}

\begin{prop}\label{c-bounded}
Let $X\subseteq\R^{d_1}$, $Y\subseteq\R^{d_2}$ be smooth submanifolds. Let $u\in(\Gen(X))^{d_2}$.
Then the following are equivalent:
\begin{enumerate}
\item $u$ is c-bounded into $Y$
\item for one, and thus for all representatives $(u_\eps)_\eps$ of $u$,
\[(\forall K\csub X)(\exists K'\csub Y)(\forall m\in\N)\big(\sup_{x\in K} d(u_\eps(x), K') = \Ord(\eps^m), \eps\to 0\big)\]
(here $d$ denotes the Euclidean distance in $\R^{d_2}$).
\item as a pointwise function on compactly generalized points, $u(\widetilde X_c)\subseteq \widetilde Y_c$.
\end{enumerate}
\end{prop}
\begin{proof}
$(1)\implies(3)$: let $\tilde x\in\widetilde X_c$. For a representative $(x_\eps)_\eps$ of $x$, $x_\eps\subseteq K\csub X$, for sufficiently small $\eps$. If $(u_\eps)_\eps$ is a representative of $u$ with $u_\eps(K)\subseteq K'\csub Y$ for sufficiently small $\eps$, then $u_\eps(x_\eps)\in K'$ for sufficiently small $\eps$, so $u(\widetilde x)\in\widetilde Y_c$.\\
$(3)\implies(2)$: suppose that there exists $K\csub X$ such that
\[
(\forall K'\csub Y)(\exists m\in\N)(\forall\eta\in (0,1))(\exists \eps < \eta)(\exists x\in K)(d(u_\eps(x), K')\ge \eps^m).
\]
We distinguish 2 cases.\\
(a) $(\sup_{x\in K} d(u_\eps(x),Y))_\eps$ is not negligible, i.e., the previous formula also holds for $Y$ itself instead of $K'$. Then we find a decreasing sequence $(\eps_n)_{n\in\N}$, with $\eps_n\to 0$ and $x_{\eps_n}\in K$, with $d(u_{\eps_n}(x_{\eps_n}), Y)\ge \eps^m$, for some $m$. Extend $(x_{\eps_n})_{n\in\N}$ to $(x_\eps)_\eps$, with $x_\eps\in K$, $\forall \eps$. Then it represents $\tilde x\in \widetilde X_c$ for which $u(\tilde x)\notin \widetilde Y_c$.\\
(b) $(\sup_{x\in K} d(u_\eps(x),Y))_\eps$ is negligible. 
Consider a compact exhaustion $(K_n)_{n\in\N}$ of $Y$ with $K_n\subseteq\interior{(K_{n-1})}$, $\forall n\in\N$.
Then we find a decreasing sequence $(\eps_n)_{n\in\N}$, with $\eps_n\to 0$, $m_n\in\N$ and $x_{\eps_n}\in K$ such that $d(u_{\eps_n}(x_{\eps_n}), Y)<\eps_n^{m_n+n}<\eps_n^{m_n}\le d(u_{\eps_n}(x_{\eps_n}), K_n)$; in particular, there exists $y_n\in Y\setminus K_n$ such that
$\abs{y_n-u_{\eps_n}(x_{\eps_n})}\le \eps_n^n$. Let $m<n$. As $K_m\subseteq\interior{(K_n)}$, $d(y_n, K_m)\ge r\in\R^+$, so $u_{\eps_n}(x_{\eps_n})\notin K_m$ as soon as $n$ is large enough.
Extend $(x_{\eps_n})_{n\in\N}$ to $(x_\eps)_\eps$, with $x_\eps\in K$, $\forall\eps$. Then it represents $\tilde x\in \widetilde X_c$ for which $u(\tilde x)\notin \widetilde Y_c$.\\
$(2)\implies(1)$: let $(u_\eps)_\eps$ be a representative of $u$. Let $W$ be a normal tubular neighbourhood of $Y$ in $\R^{d_2}$ with associated smooth retraction $q$: $W\to Y$ (see \cite{Hirsch}). By assumption, $u$ is c-bounded into $W$, so the composition $q\comp u$ is a well-defined element of $(\Gen(X))^{d_2}$ and is c-bounded into $Y$.
Let $K\csub X$. By the fact that $W$ is a normal tubular neighbourhood of $Y$, $q(x)$ is the unique element of $Y$ that is closest to $x$, for each $x\in W$. So for sufficiently small $\eps$, $\sup_{x\in K}\abs{(q\comp u_\eps)(x)-u_\eps(x)}=\sup_{x\in K}d(u_\eps(x),Y)$ which is negligible by assumption. It follows that $q\comp u = u$ as a generalized function in $(\Gen(X))^{d_2}$.
\end{proof}

\begin{cor}\label{cor c-bounded}
$(1)$ Let $X\subseteq\R^{d_1}$, $Y\subseteq\R^{d_2}$ be smooth submanifolds.
An element $u\in(\Gen(X))^{d_2}$ that is c-bounded into $Y\subseteq\R^{d_2}$ defines a unique element of $\Genld[X,Y]$ by restricting a representative satisfying eqn.~(\ref{eqn c-bounded}) to (suitably chosen, depending on $\eps$) compact subsets of $X$.\\
$(2)$ Let $X$, $Y$ be smooth manifolds. Let $u\in\Genld[X,Y]$ and $v\in\Gen(Y)$. Then the composition $v\comp u$ defined on representatives by means of
$(v\comp u)_\eps = v_\eps\comp u_\eps$
is a well-defined generalized function in $\Gen(X)$.
\end{cor}
\begin{proof}
(1) Let $(u_\eps)_\eps$ be a representative of $u$ satisfying eqn.~(\ref{eqn c-bounded}). Let $(K_n)_{n\in\N}$ be a compact exhaustion of $X$. Then for each $n\in\N$, there exists $K\csub Y$ and $\eps_n\in (0,1)$ such that $u_\eps(K_n)\subseteq K$, for each $\eps\le\eps_n$. We may suppose $(\eps_n)_{n\in\N}$ to be a decreasing sequence. Let $v_\eps=\restr{u_\eps}{K_n}$, for each $\eps_{n+1}<\eps\le\eps_n$. Then $(v_\eps)_\eps$ represents an element of $\Genld[X,Y]$. Well-definedness 
follows as in \cite[Prop.~3.2.43]{GKOS}.\\
(2) Analogous to \cite[Prop.~3.2.58]{GKOS}.
\end{proof}
It follows that, in case $\Gen[X,Y]\subsetneq\Genld[X,Y]$, a characterization of algebra homomorphisms $\Gen(X)\to\Gen(Y)$ as compositions with generalized maps $\Gen[X,Y]$ is not possible ($\Genld[X,Y]$ has to be used instead).

\end{document}